# Subgroups of Symmetric Groups Determined Using Signed Brauer Diagrams

Ram Parkash Sharma and Rajni Parmar

Department of Mathematics, Himachal Pradesh University
Summerhill, Shimla (India)
Email: rp_math_hpu@yahoo.com

**Abstract**

The group $\overrightarrow{S_n}$ of signed Brauer diagrams with no horizontal edges becomes a subgroup of $S_{2n}$. Therefore, it is natural to find the elements of $\overrightarrow{S_n}$ which are even permutations in $S_{2n}$. There are $2^{n-1}n!$ such elements that form a subgroup of $S_{2n}$. Consequently, $S_n$ has a class of subgroups of order 1, $2^{2-1}2!$, $2^{3-1}3!, \ldots,$ $2^{\left(\left[\frac{n}{2}\right]-1\right)}\left[\frac{n}{2}\right]!$ of even permutations and hence $2^{\left[\frac{n}{2}\right]}\left[\frac{n}{2}\right]! \mid n!$. Finally, using these results it is shown that every group of order $n(n-1)\ldots(\left[\frac{n}{2}\right]+1)$ has a $Sylow\ 2-\ subgroup$ of order $2^{\left[\frac{n}{2}\right]}$.

## 1. Introduction

Parvathi and Kamaraj [2] defined signed Brauer algebras $\overrightarrow{B_n}(x)$ over $\mathbb{C}$, $(x \in \mathbb{C})$ having basis consisting of signed diagrams $\vec{d}$, which consist of two rows of $n$ points labelled $\{1,2,\ldots,n\}$ in first row and $\{(n+1),(n+2),\ldots,2n\}$ in second row, with each dot joined to precisely one other dot (distinct from itself). Each edge of the diagram $\vec{d}$ has a direction: either positive or negative. A positive vertical (horizontal) edge is denoted by $\downarrow$ $(\rightarrow)$ and a negative vertical (horizontal) edge is denoted by $\uparrow$ $(\leftarrow)$. An element in $\overrightarrow{B_5}(x)$ is given as:

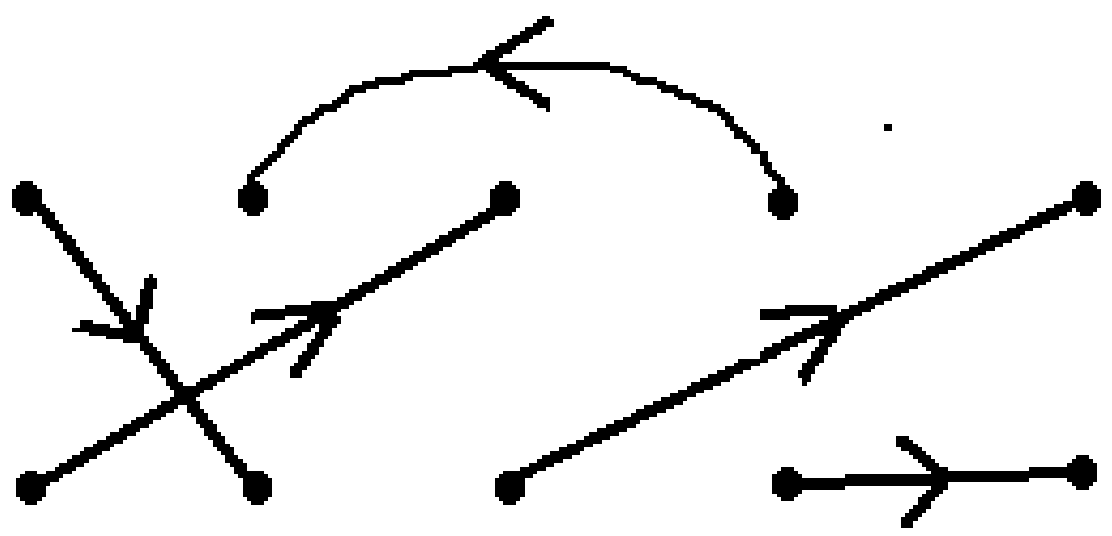

To define multiplication in $\overrightarrow{B_n}(x)$, it is enough to define the product of two signed diagrams $\overrightarrow{d_1}$ and $\overrightarrow{d_2}$, to do so place $\overrightarrow{d_1}$ above $\overrightarrow{d_2}$ and join corresponding points (interior signed loops are deleted). Let $\overrightarrow{d_3}$ be the signed diagram obtained in this way. An edge in $\overrightarrow{d_3}$ is positive (negative) according as the number of negative edges involved to form this edge is even (odd). Similarly, a loop β in $\overrightarrow{d_1} o \overrightarrow{d_2}$ is so called positive (negative).

The multiplication is then defined

$$\overrightarrow{d_1} o \overrightarrow{d_2} = x^{2n_1+n_2}\overrightarrow{d_3}$$

where $n_1$ and $n_2$ are the numbers of positive and negative loops respectively formed placing $\overrightarrow{d_1}$ above $\overrightarrow{d_2}$.

The set of signed Brauer diagrams $\overrightarrow{B_n}$ contains the set of Brauer diagrams $B_n$ (each Brauer diagram corresponds to same underlying diagram with all positive edges). For more detail on Brauer and signed Brauer algebras, one can refer to [3,4,5,6]. The Brauer diagrams with no horizontal edges form a permutation group, by identifying a Brauer diagram with a permutation $\pi$ connecting the $i^{th}$ upper vertex to the $\pi(i)^{th}$ lower vertex. In a similar way, the signed Brauer diagrams with no horizontal edges form a group $\overrightarrow{S_n}$ which corresponds to the wreath product of $Z_2$ by $S_n$ [see 1]. The group $\overrightarrow{S_n}$ becomes a subgroup of $S_{2n}$. In order to identify $\overrightarrow{S_n}$ as a subgroup of $S_{2n}$, first divide the set

$$\{1,2,\ldots n,\ldots,2n\}$$

into two equal parts that is the numbers $\leq n$ and $> n$. To distinguish the numbers less than or equal to $n$ and greater than $n$, we fix the notation $r_1, r_2, \ldots, r_n$ for the numbers less than or equal to $n$ and $s_1, s_2, \ldots, s_n$ for the numbers greater than $n$. We fix another notation as:

For each $i \in \{1,2,\ldots,2n\}$, we denote by $i^*$ the number given by:

$$i^* = \begin{cases} i+n, & 1 \le i \le n, \\ i-n, & n < i \le 2n \end{cases}.$$

Thus $r_i^* > n$; $i = 1,2,\ldots,n$ and $s_j^* \le n$; $j = 1,2,\ldots,n$.

The elements of $\overrightarrow{S_n}$ are identified in $S_{2n}$ as follows:

Let $\sigma \in \overrightarrow{S_n}$ have a negative edge

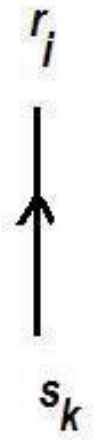

between $r_i$ and $s_k$. Then the number $r_i(\le n) \longmapsto s_k(> n)$ and $r_i^* \longmapsto s_k^*$ under $\sigma$ in $S_{2n}$ and vice-versa. If we have a positive edge

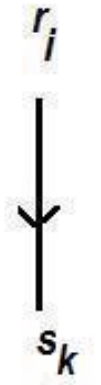

between $r_i$ and $s_k$. Then $r_i \longmapsto s_k^*$ and $r_i^* \longmapsto s_k$.

**Example 1.1.** Consider the eight elements of $\overrightarrow{S_2}$ given below:

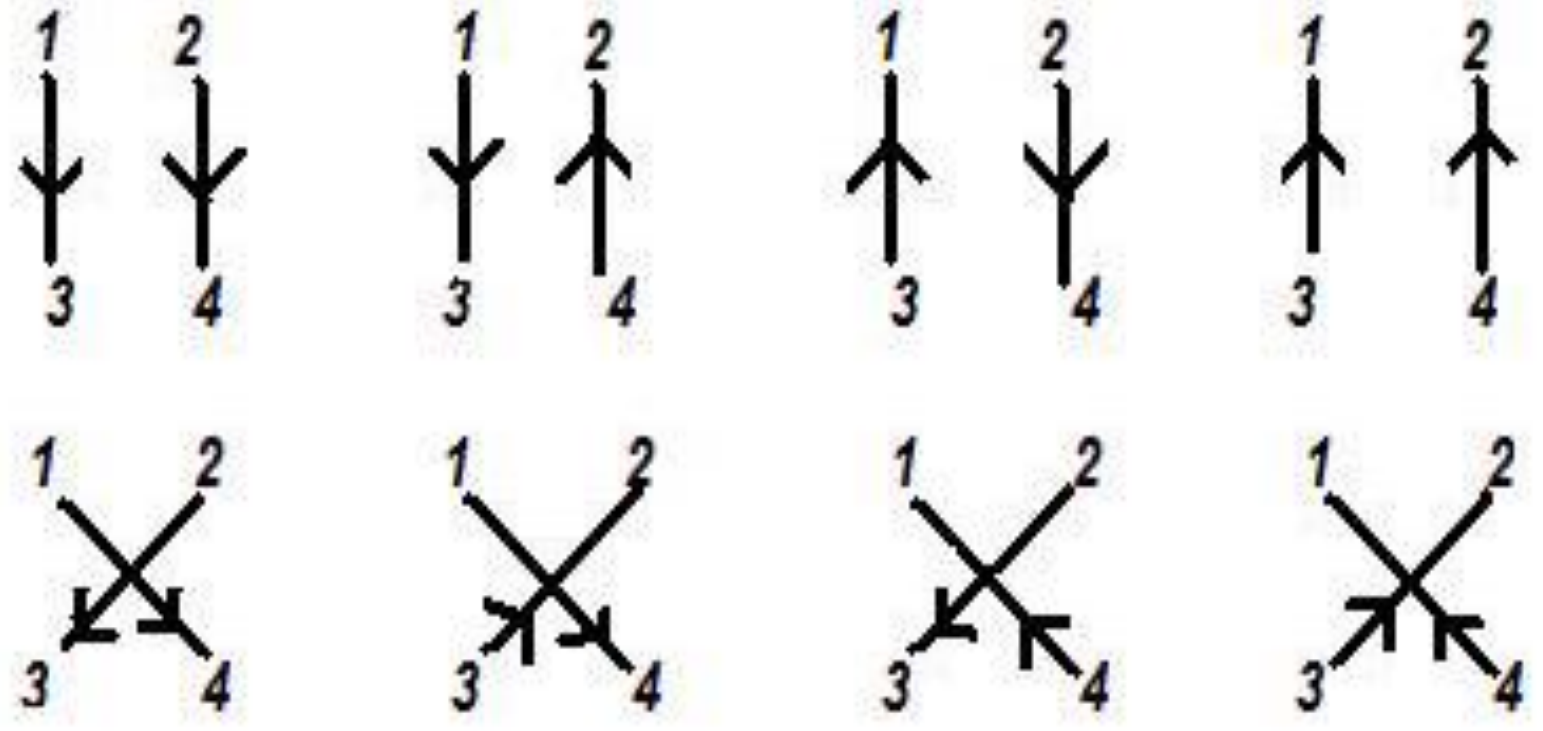

According to the above definition, the elements of $\overrightarrow{S_2}$ are identified in $S_4$ as follows:

$$\{e, (2\ 4), (1\ 3), (1\ 3)(2\ 4), (1\ 2)(3\ 4), (1\ 2\ 3\ 4), (1\ 4\ 3\ 2), (1\ 4\ )(2\ 3)\}.$$

As $\overrightarrow{S_n}$ is a subgroup of $S_{2n}$, we started this paper with the aim of finding the elements of $\overrightarrow{S_n}$ which correspond to even permutations in $S_{2n}$ and prove that the elements of $\overrightarrow{S_n}$ having even number of negative edges (see Theorem 3.1) are even permutations in $S_{2n}$. Since the product of two negative edges in $\overrightarrow{B_n}$ gives a positive edge and product of negative and positive edges is a negative edge, it is easy to see that these elements form a subgroup of $\overrightarrow{S_n}$ and hence of $S_{2n}$. The total number of elements with even number of negative directions in $\overrightarrow{S_n}$ is given by:

$$\begin{cases} ({}^n c_0 + {}^n c_2 + {}^n c_4 + {}^n c_6 + \dots + {}^n c_n) n!, \text{if } n \text{ is even} \\ ({}^n c_0 + {}^n c_2 + {}^n c_4 + \dots + {}^n c_{n-1}) n!, \text{if } n \text{ is odd} \end{cases}$$
$$= 2^{n-1} n!.$$

Thus, $S_{2n}$ has a subgroup of order $2^{n-1}n!$. In particular, $S_n$ has a subgroup of order $2^{\left(\left[\frac{n}{2}\right]-1\right)}\left[\frac{n}{2}\right]!$ for $n \geq 2$ and hence in general it has a class of subgroups of order 1, $2^{2-1}2!$, $2^{3-1}3!, \dots, 2^{\left(\left[\frac{n}{2}\right]-1\right)}\left[\frac{n}{2}\right]!$ of even permutations. In consequence of it, $2^{\left[\frac{n}{2}\right]}\left[\frac{n}{2}\right]! \mid n!$; hence

$$2^{\left[\frac{n}{2}\right]} \mid n(n-1)\dots\left(\left[\frac{n}{2}\right]+1\right).$$

Further, we observe that (Corollary 3.4)

$$2^{\left(\left[\frac{n}{2}\right]+1\right)} \nmid n(n-1)\dots\left(\left[\frac{n}{2}\right]+1\right).$$

Therefore, every group of order $n(n-1)\dots(\left[\frac{n}{2}\right]+1)$ has a $Sylow\ 2-\ subgroup$ of order $2^{\left[\frac{n}{2}\right]}$.

## 2. Even Permutations of $\overrightarrow{S_n}$ in $S_{2n}$

First, we prove that the elements of $\overrightarrow{S_n}$ having positive direction of all the edges are even permutations in $S_{2n}$. We, very often, use $r_i \mapsto s_j$ for $\sigma(r_i) = s_j$ and $r_i \not\mapsto s_j$ for $\sigma(r_i) \neq s_j$ when $\sigma$ is understood.

**Lemma 2.1**. If $\sigma \in \overrightarrow{S_n}$ has positive direction of all the edges, then $\sigma$ is an even permutation in $S_{2n}$.

**Proof**. Let $\sigma \in \overrightarrow{S_n}$ has all edges of positive direction. Then $r_i \leq n$, goes to some $r_j = (s_j^*) \leq n$ under $\sigma$. Consequently, $r_i^*$ goes to some $r_j^* = (s_j^*)^* = s_j$ and both $r_i^*$,

$s_j > n.$ In this way, we get

$$r_1 \mapsto r_2 \mapsto r_3 ... \mapsto r_i,$$

and

$$r_1^* \mapsto r_2^* \mapsto r_3^* ... \mapsto r_i^*.$$

Note that $r_i \not\mapsto r_1^*$ and $r_i^* \not\mapsto r_1$ because this is only possible if we have a negative edge between $r_i$ and $r_1^*$. Since direction of all the edges of $\sigma$ is positive, there is only one possibility of having $r_i \mapsto r_1$ and hence $r_i^* \mapsto r_1^*$ using a positive edge between $r_i$ and $r_1^*$. Therefore, corresponding to every cycle $(r_1 \ r_2 \ r_3 ... r_i)$ written using the numbers $\leq n,$ there is another cycle $(r_1^* \ r_2^* \ r_3^* ... r_i^*)$ of the numbers greater than $n$ having same length. So $\sigma$ is an even permutation in $S_{2n}$.

**Example 2.2.** Consider an element

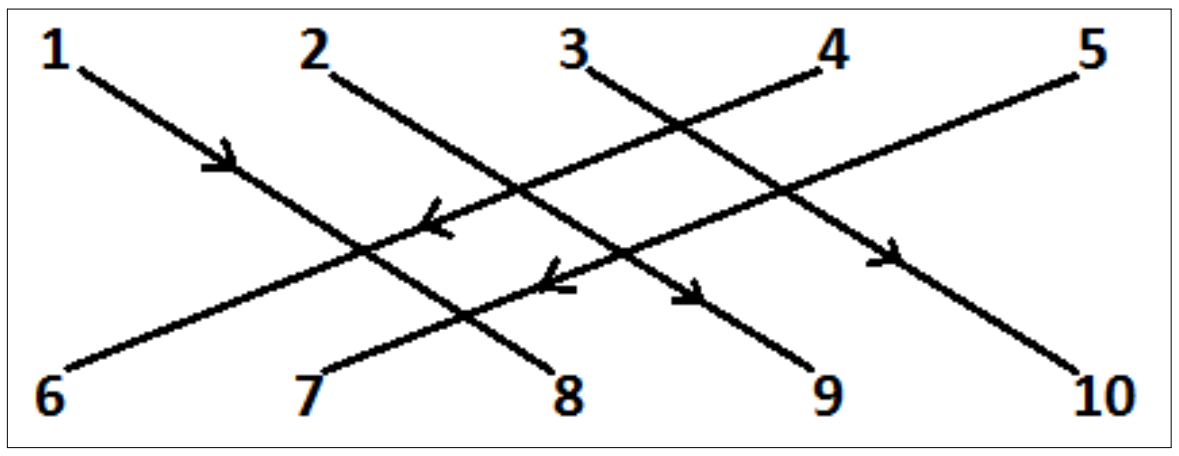

of $\overrightarrow{S_5}$ having all edges of positive direction. This element becomes

$$(1 \ 3 \ 5 \ 2 \ 4)(6 \ 8 \ 10 \ 7 \ 9)$$

in $S_{10}$.

For $\sigma \in \overrightarrow{S_n}$ with negative direction of all the edges, any cycle $c_i$ , in the product of cycles of $\sigma,$ no two numbers $r_i$ , $r_j$ occur together and same is true for $s_i$ , $s_j$ also. With this observation, we have

**Lemma 2.3.** Let $\sigma \in \overrightarrow{S_n}$ has negative direction of all the edges. Then, length of each cycle in the product of $\sigma$ in $S_{2n}$ is always even.

**Proof.** Since all the edges are of negative direction, any number $r_i \leq n$ goes to a number $s_j > n.$ Since $s_j^* \leq n,$ $s_j^*$ can be $r_i$. In this case, we get a transposition $(r_i \, r_i^*)$. Let $(r_1 \ s_1 ...)$ be a cycle in the product of $\sigma$. In this cycle no $r_i, r_j$ or $s_i, s_j$ occur together as all the edges are of negative direction. Moreover the last number will be some $s_k > n.$ For, if we have a cycle of the type $(r_1 \ s_1 \ r_2 \ s_2 ... \ r_k)$ . Then $r_k \mapsto r_1$ gives a positive edge between $r_k$ and $r_1^*$. Hence any cycle in the product of $\sigma$ will be of the type $(r_1\ s_1\ r_2\ s_2 \ldots r_{k-1}\ s_k)$ having length $2k$.

The above result does not depend on whether $n$ is even or odd. However, the number of cycles in the product of $\sigma \in \overrightarrow{S_n}$ having negative direction of all the edges depends upon $n$ as observed in the following examples.

**Example 2.4.** For $n$ to be odd, consider an element of $\overrightarrow{S_5}$ having all edges of negative direction. This element becomes $(1\ \ 7\ \ 3\ \ 9)(6\ \ 2\ \ 8\ \ 4)(5\ \ 10)$ in $S_{10}$. Clearly, length of each cycle in the product of $\sigma$ is even, but number of cycles is 3.

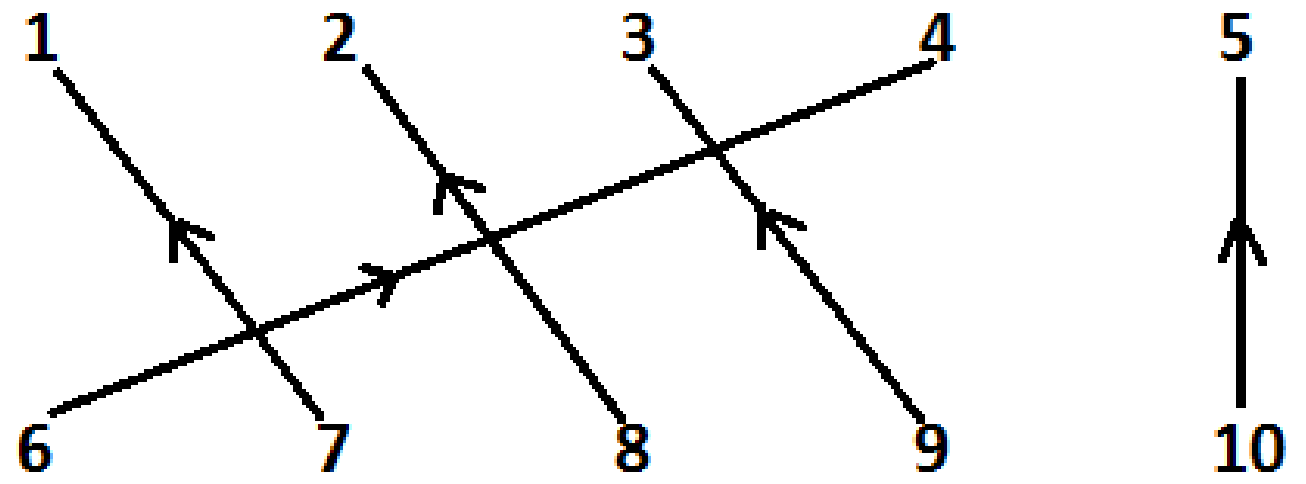

**Example 2.5**. For $n$ to be even, consider an element

$\sigma =$

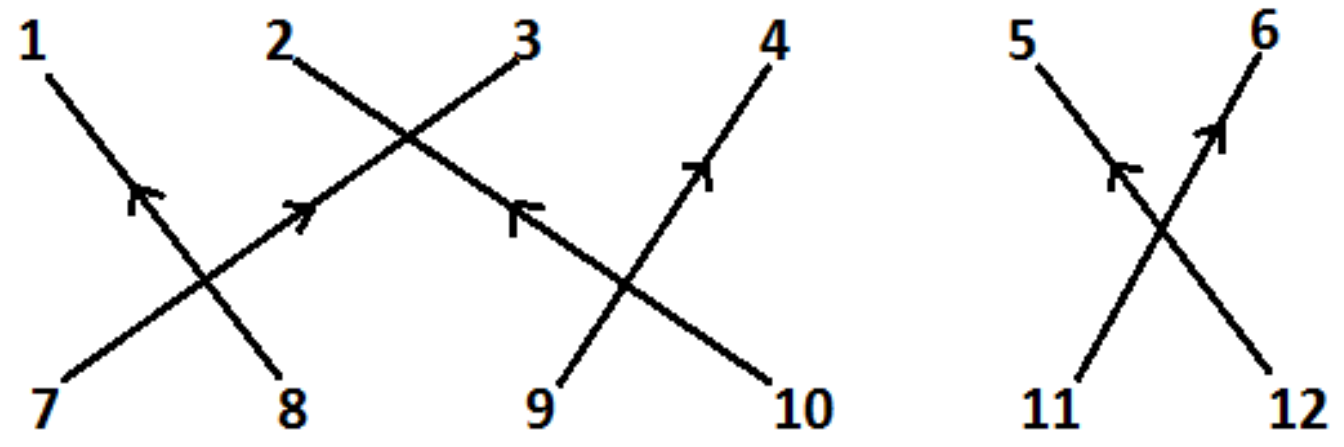

of $\overrightarrow{S_6}$ having all edges of negative direction which becomes $(1\ \ 8\ \ 4\ \ 9)(7\ \ 2\ \ 10\ \ 3)(5\ \ 12)(6\ \ 11)$ in $S_{12}$. Here, again the length of each cycle is even, but the number of cycles is also even.

With these observations, we prove a general result regarding the number of cycles of $\sigma \in \overrightarrow{S_n}$ having negative direction of all the edges.

**Lemma 2.6.** Let $\sigma \in \overrightarrow{S_n}$ has negative direction of all the edges. Then the number of disjoint cycles in the product of $\sigma$ is even if $n$ is even and odd if $n$ is odd.

**Proof.** For given $\sigma \in \overrightarrow{S_n}$ with negative direction of all the edges, note that $\sigma$

fixes no $r_i$ or $s_j$, as $\sigma(r_i) = r_i$ if and only if we have a positive edge between $r_i$ and $r_i^*$. Similarly $\sigma(s_j) = s_j$ if and only if we have a positive edge between $s_j^*$ and $s_j$. Therefore all the numbers $1, 2, 3, ..., n, n+1, ...2n$ will be involved in the cycles of $\sigma$. We discuss different cases for involvement of even or odd number of numbers in a cycle of $\sigma$. Here we also use the fact observed in Lemma 2.3 that in any cycle of $\sigma \in \overrightarrow{S_n}$, no $r_i$, $r_j$ or $s_i$, $s_j$ occur together as all the edges are of negative direction.

**Case (i)** Suppose we want to involve two numbers from each side of $n$ in the product of cycles of $\sigma$, then the only possibility is $r_i \mapsto s_j$ $(s_j \neq r_i^*)$ and $r_i^* \mapsto s_j^*$ . In this case, we get two cycles $(r_i\ s_j)$ and $(r_i^*\ s_j^*)$.

The possibility of having these four numbers in a single cycle is ruled out because these four numbers will be in one cycle if and only if we have a positive edge between $s_j^*$ and $r_j^*$.

Further, if we want to use four numbers from each side of $n$, then

$$r_1 \mapsto s_1 \mapsto r_2 \mapsto s_2,$$

and

$$r_1^* \mapsto s_1^* \mapsto r_2^* \mapsto s_2^*.$$

(Note that without loss of generality, we can start with $r_1$, because if we start with $s_1$ then we are writing the second row in place of first). Here we take $s_i \not\mapsto r_i$, otherwise we get two cycles of length 2, which we have already discussed). These sequences will terminate if $s_2 \mapsto r_1$ and $s_2^* \mapsto r_1^*$, because $s_2 \not\mapsto r_1^*$, $s_2^* \not\mapsto r_1$ being on the same side of $n$ as $\sigma$ has all negative edges. So we get two cycles

$$(r_1 s_1 r_2\ s_2) \text{ and } (r_1^*\ s_1^* r_2^* s_2^*).$$

In general, the involvement of an even number of numbers from each side is possible if and only if we take

$$r_1 \mapsto s_1 \mapsto r_2 \mapsto s_2 \mapsto ... \mapsto r_i \mapsto s_i \quad ,$$

and

$$r_1^* \mapsto s_1^* \mapsto r_2^* \mapsto s_2^* \mapsto ... \mapsto r_i^* \mapsto s_i^*.$$

Again, in this case, we get two cycles $(r_1\ s_1\ r_2\ s_2\ \ldots .r_i\ s_i)$ and $(r_1^*\ s_1^*$ $r_2^*\ s_2^*.\ldots .r_i^*\ s_i^*)$.

**Case (ii)** Now suppose that only one number from either side of $n$ is involved in a cycle of $\sigma$ , then $r_i \mapsto r_i^*$ and we get only one cycle $(r_i \quad r_i^*)$ , a negative transposition.

Three numbers can be involved from either side of $n$ in a cycle of $\sigma$ if and only if

$$r_1 \mapsto s_1 \mapsto r_2, (s_1 \neq r_1^* \text{ and } r_2 \neq s_1^*)$$

and

$$r_1^* \mapsto s_1^* \mapsto r_2^*.$$

These sequences will terminate if $r_2 \mapsto r_1^*$ and $r_2^* \mapsto r_1$, because $r_2 \not\mapsto r_1$ and $r_2^* \not\mapsto r_1^*$ being on the same side of $n$. So, we get a single cycle $(r_1\, s_1\ r_2$ $r_1^*\ s_1^*\, r_2^*)$.

In general, if an odd number of $r_i$ and $s_i$ are to be involved from either side of $n$ in a cycle of $\sigma$, then

$$r_1 \mapsto s_1 \mapsto r_2 \mapsto s_2 \mapsto r_3 \mapsto s_3 \mapsto r_4.\ldots .\mapsto s_{i-1} \mapsto r_i$$

and

$$r_1^* \mapsto s_1^* \mapsto r_2^* \mapsto s_2^* \mapsto r_3^* \mapsto s_3^* \mapsto r_4 \ldots \mapsto s_{i-1}^* \mapsto r_i^*.$$

These sequences will terminate if $r_i \mapsto r_1^*$ and $r_i^* \mapsto r_1$, because $r_i \not\mapsto r_1$ and $r_i^* \not\mapsto r_1^*$ being on the same side of $n$. So we get a single cycle $(r_1\, s_1 r_2$ $s_2 \ldots s_{i-1} r_i\ r_1^*\ s_1^*\ r_2^*\, s_2^* \ldots .s_{i-1}^*\, r_i^*)$ of even length.

Therefore we conclude that there are two type of cycles in the product of $\sigma$ :

(i) cycles having even number of $r_i \leq n$ and

(ii) cycles having odd number of $r_i \leq n$.

The cycles of first kind exist in pair. A cycle of second kind, having length $2t$, corresponds to some odd number $t$ which exists in $\lambda \vdash n$. Therefore if $n$ is even, then corresponding to an odd number $t$ in the partition $\lambda \vdash n$ , there exists $t'$ an odd number in the same partition $\lambda \vdash n$. Hence we get a pair of cycles of type $2t$ and $2t'$ (see case (ii)). Hence, the number of cycles in the product of $\sigma$ will be even if $n$ is even and odd if $n$ is odd.

As an easy consequence of the above result, we have

**Corollary 2.7.** Any cycle in the product of $\sigma \in \overrightarrow{S_n}$ having all the edges of negative direction, ends with some $s_j$ if it starts with $r_i$.

**Corollary 2.8.** Let $\sigma \in \overrightarrow{S_n}$ with $n$ even, has negative direction of all the edges. Then, $\sigma$ is not a $2n$ –cycle in $S_{2n}$.
**Proof.** This falls in the case (i) of the above theorem. In this case, for each cycle, there corresponds another cycle of the same length, in the product of $\sigma$. Hence any such $\sigma$ cannot be a single cycle of length $2n$.

Using the above results proved for $\sigma \in \overrightarrow{S_n}$ having negative direction of all the edges, we have
**Proposition 2.9.** Let $\sigma \in \overrightarrow{S_n}$ with $n$ even, has negative direction of all the edges. Then $\sigma$ is an even permutation in $S_{2n}$.
**Proof.** By lemma 2.6, number of disjoint cycles of $\sigma$ in $S_{2n}$ is even if and only if $n$ is even. Also by lemma 2.3, length of each cycle of $\sigma$ is even. Even number of disjoint cycles with even length will give even number of transpositions, hence $\sigma$ is an even permutation in $S_{2n}$.

## 3. Main Results

Now we prove the main results of this paper.
**Theorem 3.1.** For $n \geq 2$, let $n = 2m + r, 0 \leq m \leq \left[\frac{n}{2}\right]$ and $0 \leq r \leq n$. Then $\sigma \in \overrightarrow{S_n}$ is an even permutation in $S_{2n}$ if and only if $\sigma$ has $2m$ edges of negative direction.
**Proof.** First, we analyze the role of positive and negative edges to get a cycle in the product of $\sigma$ in $S_{2n}$.
One negative edge between $r_i$ and $r_i^*$ gives a transposition $(r_i \;\; r_i^*)$ in $S_{2n}$. Thus, if we have such $2m-$ negative edges, then we get an even number of transpositions $(r_1 \;\; r_1^*), (r_2 \;\; r_2^*), (r_3 \;\; r_3^*), ..., (r_{2m} \;\; r_{2m}^*)$ in $S_{2n}$ corresponding to these negative edges. The remaining positive edges give pair of cycles of equal length as observed in Lemma 2.1 and hence $\sigma$ is an even permutation in $S_{2n}$.
Suppose some mixed positive and negative edges give a cycle in the product of $\sigma$. Note that two numbers $r_i, \; r_j$ or $s_i, \; s_j$ in a cycle occur if and only if there is a positive edge and two numbers $r_i, \; s_j$ or $s_i, \; r_j$ occur if and only if there is a negative edge.
**Case (i)** Suppose a sequence of the type

$$r_1 \mapsto .... \mapsto s_i ... \mapsto r_j$$

occurs inside a cycle of $\sigma$. In this sequence, there are even number of changes from $r_i$ to $s_i$ and $s_j$ to $r_i$ (possibly no change). As observed above this sequence occurs if and only if there exists a sequence of even number of negative edges (possibly no negative edge). But then we also have

$$r_1^* \mapsto \ldots \mapsto s_i^* \ldots \mapsto r_j^*.$$

In order to complete the cycles either $r_j \mapsto r_1$ and hence $r_j^* \mapsto r_1^*$ or $r_j \mapsto r_1^*$ and hence $r_j^* \mapsto r_1$. Here we are not taking the case $r_j \mapsto r_j^*$ because that will remove $r_j$ and $r_j^*$ from this cycle and single cycles have been discussed earlier.

**Subcase (i)** $r_j \mapsto r_1$ and $r_j^* \mapsto r_1^*$. This is only possible with an additional positive edge between $r_j$ and $r_1^*$. Therefore the number of negative edges used remain even. Hence in this case, we have used even number of negative edges with some positive edges and we get two cycles $(r_1 \ldots . s_i \ldots \ r_j)$ and $(r_1^* \ldots . s_i^* \ldots . r_j^*)$ contributing an even number of transpositions in the product of $\sigma$.

**Subcase (ii)** $r_j \mapsto r_1^*$ and $r_j^* \mapsto r_1$. This is only possible with an additional negative edge between $r_j$ and $r_1^*$. Hence in this case, number of negative edges used changes to odd and we get a single cycle $(r_1 \ldots . s_i \ldots \ r_j \ r_1^* \ldots . s_i^* \ldots . r_j^*)$. Clearly the length of this cycle is even as number of $r_i^*$ = number of $r_i$ and number of $s_i^*$ = number of $s_i$. This contributes an odd number of transpositions in the product of $\sigma$.

**Case (ii)** Now consider a sequence of the type

$$r_1 \mapsto \ldots \mapsto r_i \ldots \ \mapsto \ \ s_i \ \mapsto \ldots \mapsto s_j.$$

Here, there are odd number of changes from $r_i$ to $s_i$ and $s_j$ to $r_j$ (both included). Hence this sequence occurs inside a cycle of $\sigma$ if and only if there exists a sequence of odd number of negative edges. But then we also have

$$r_1^* \mapsto \ldots \mapsto r_i^* \ldots \ \mapsto s_i^* \ \mapsto \ldots \mapsto s_j^*.$$

Again in order to complete the cycle either $s_j \mapsto r_1$ and hence $s_j^* \mapsto r_1^*$ or $s_j \mapsto r_1^*$ and hence $s_j^* \mapsto r_1$.

**Subcase (i)** $s_j \mapsto r_1$ and $s_j^* \mapsto r_1^*$. This is only possible with an additional negative edge between $s_j^*$ and $r_1^*$. Hence, in this case, number of negative edges used, changes to even and we get two cycles $(r_1 \ldots . r_i \ldots . s_i \ \ldots . s_j)$ and $(r_1^* \ldots . r_i^* \ldots \ s_i^* \ \ldots s_j^*)$ contributing an even number of transpositions in the product of $\sigma$.

**Subcase (ii)** $s_j \mapsto r_1^*$ and $s_j^* \mapsto r_1$. This is possible only with an additional positive edge between $s_j^* \mapsto r_1^*$. In this case number of negative edges used remains odd and we get a cycle

$$(r_1 \ldots . r_i \ldots . s_i \ \ldots . s_j \ \ r_1^* \ldots . r_i^* \ldots \ s_i^* \ \ldots . s_j^*)$$

of even length. This contributes in product of $\sigma$ an odd number of transpositions.

**Conclusion:** With the above observations, finally we conclude that an even number of negative edges (possibly zero) together with some positive edges contribute two cycles of same length in the product of $\sigma$ and an odd number of negative edges (possibly one) contribute a product of odd number of transpositions for $\sigma$ . If we take any odd number $r$ in a partition of $2m$ , then there exists an another odd number $r'$ in this partition. Hence, by above discussion it follows that in the product of $\sigma,$ there are pairs of cycles of equal length contributing to even number of transposition and some more even number of transpositions coming from the pairs $(r , r')$ of negative edges ( $r$ and $r'$ are odd numbers). Therefore, $\sigma$ is an even permutation in $S_{2n}$.

Let $\sigma \in \overrightarrow{S_n}$ has $2k+1$ $(k \geq 1)$ edges with negative direction. Then in any partition of $2k+1$, the number of odd numbers will be odd. Therefore, from the above discussion, the total number of transposition in the product of $\sigma$ in $S_{2n}$ will odd.

**Theorem 3.2.** Let $\overrightarrow{A_n} = \{\sigma \in \overrightarrow{S_n} \mid \sigma \text{ has even number of negative edges}\}$. Then $\overrightarrow{A_n}$ is a normal subgroup of $\overrightarrow{S_n}$ of order $2^{n-1}n!$. Further $\overrightarrow{A_n} = \overrightarrow{S_n} \cap A_{2n}$ and hence $S_{2n}$ has a subgroup of order $2^{n-1}n!$. In particular, $S_n$ has a subgroup of order $2^{\left(\left[\frac{n}{2}\right]-1\right)}\left[\frac{n}{2}\right]$ of even permutations.

**Proof.** The multiplication defined in $\overrightarrow{B_n}$ yields that $\overrightarrow{A_n} \leq \overrightarrow{B_n}$. The result $\overrightarrow{A_n} = \overrightarrow{S_n} \cap A_{2n}$ follows from Theorem 3.1. The order of $\overrightarrow{A_n}$ is $2^{n-1}n!$ as observed in the Introduction. Also, $[\overrightarrow{S_n} : \overrightarrow{A_n}] = 2$ implies that $\overrightarrow{A_n} \lhd \overrightarrow{S_n}$ . Since $S_{2n}$ is a subgroup of $S_{2n+1,}$ we have a more general result, $S_n$ has subgroups of order 1, $2^{2-1}2!,\ 2^{3-1}3!, \ldots, 2^{\left(\left[\frac{n}{2}\right]-1\right)}\left[\frac{n}{2}\right]!$ for $n \geq 2$.

**Example 3.3.** For $k$ = 2 and 3, we have 2=2.0+2, 2=2.1+0 and 3=2.0+3, 3=2.1+1. Therefore, $S_6$ has two subgroups $\overrightarrow{A_2}$ and $\overrightarrow{A_3}$ of order $2_{C_0}2! + 2_{C_2}2! = 2^{2-1}2! = 4$ and $3_{C_0}3! + 3_{C_2}3! = 2^{\left[\frac{6}{2}\right]-1}\left[\frac{6}{2}\right]! = 2^2[3]! = 24$ respectively, which are obtained from the elements of $\overrightarrow{S_2}$ and $\overrightarrow{S_3}$ having even number of negative edges.

The subgroup $\overrightarrow{A_2}$ can be written from Example 1.1; that is,

$\{e,(13)(24),(12)(34),(14)(23)\}$. Here, we find the subgroup $\overrightarrow{A_3}$ as follows:

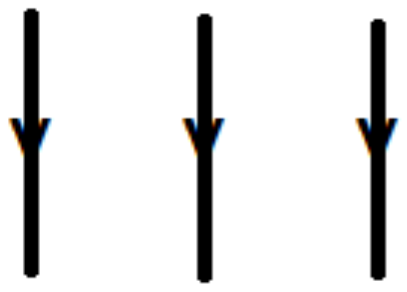
= e,

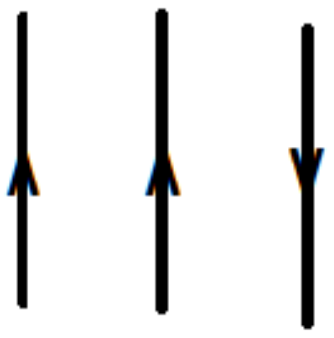
= (14)(25),

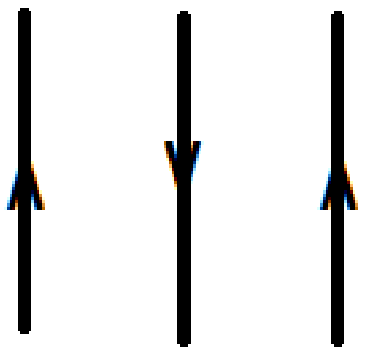
= (14)(36),

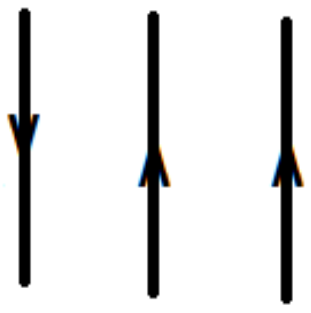
= (25)(36),

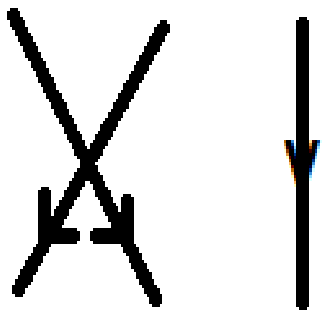
=(12)(45),

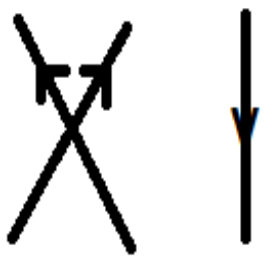
= (15)(24),

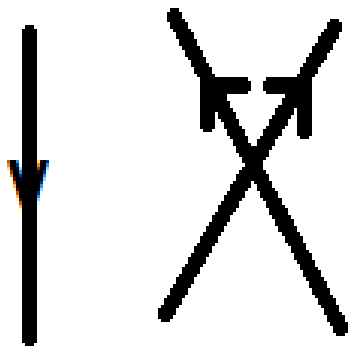
= (26)(35),

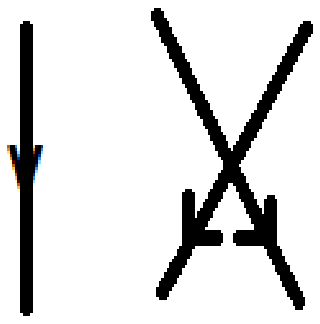
= (23)(56),

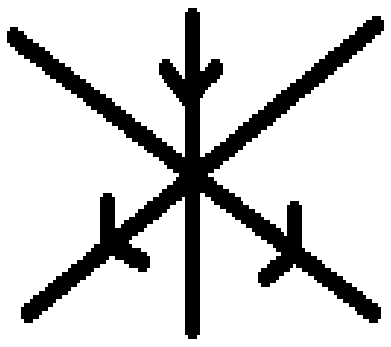
= (13)(46),

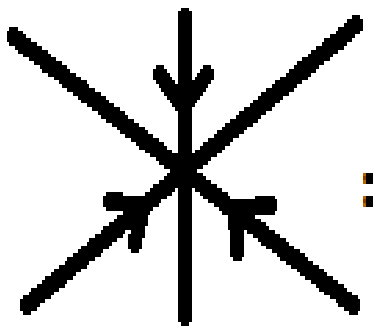
= (16)(34),

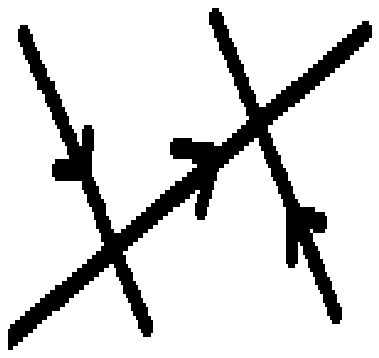
= (126)(345),

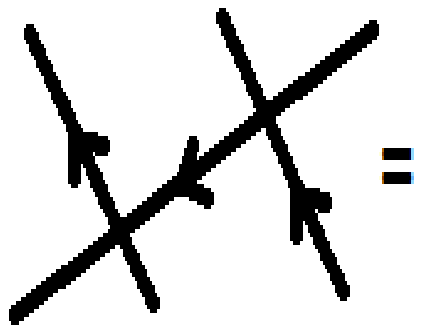
= (153)(264),

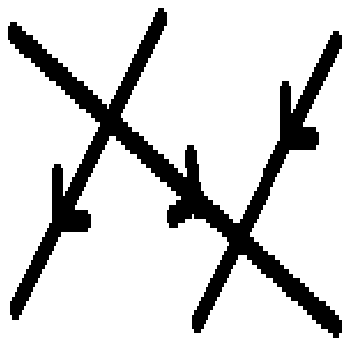
= (132)(465),

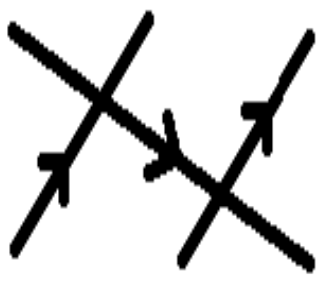
= (135)(246),

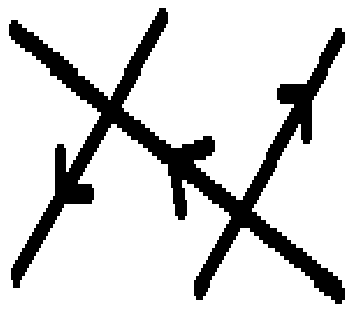 = (162)(354), 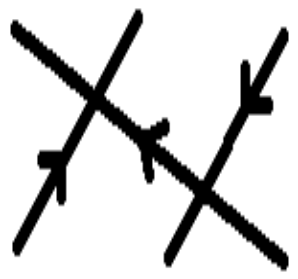 = (165)(243),

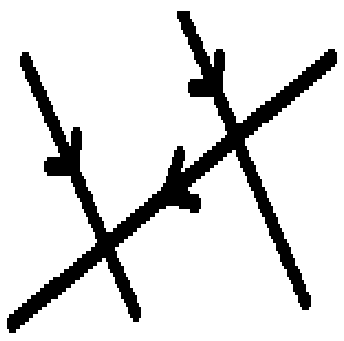 = (123)(564), 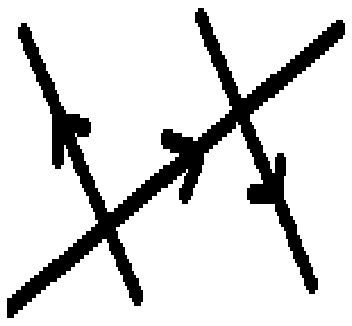 = (234)(156),

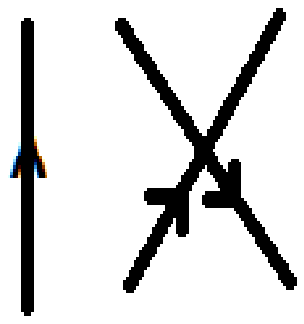 =(14)(2356), 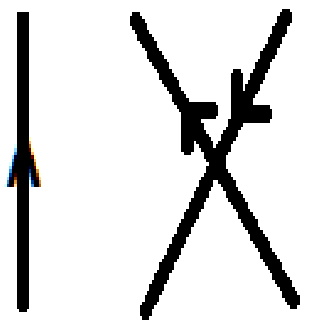 = (14)(2653),

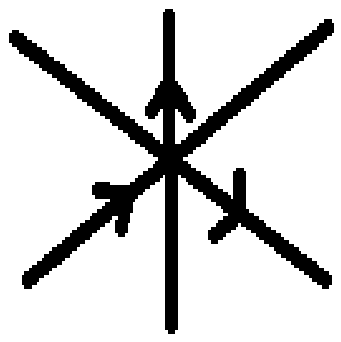 = (25)(1346), 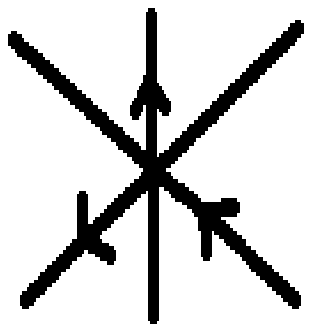 = (25)(1643),

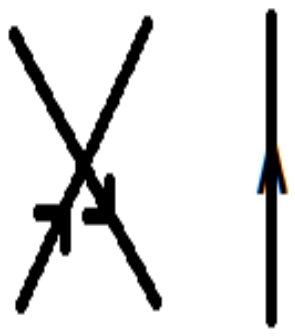 = (36)(1245), 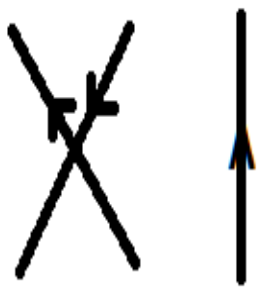 = (36)(1542)

Finally, we obtain a result as an application of Theorem 3.2.

**Corollary 3.4.** Every group of order $n(n-1)\dots([\frac{n}{2}]+1)$ has a $Sylow\ 2-subgroup$ subgroup of order $2^{[\frac{n}{2}]}$.

**Proof.** Using theorem 3.2, $2^{[\frac{n}{2}]}[\frac{n}{2}]!\ |\ n!$, hence
$2^{([\frac{n}{2}]-1)}\ |\ n(n-1)\dots\left([\frac{n}{2}]+1\right)$. We claim that

$$2^{([\frac{n}{2}]+1)} \nmid n(n-1)\dots\left([\frac{n}{2}]+1\right), n\geq 2.$$

For, we show that

$$\left([\tfrac{n}{2}]+1\right)\left([\tfrac{n}{2}]+2\right)\dots n = 2^{[\frac{n}{2}]}k, \text{ for some positive odd integer } k.$$

The result is true for $n = 2$.

Suppose that the result is true for all integers $t, 2 \leq t \leq n$.

Case(i) If $n$ is even, then

$$2^{\frac{n}{2}}k = \left(\frac{\mathrm{n}}{2}+1\right)\left(\frac{\mathrm{n}}{2}+2\right)\dots n$$

implies

$$2^{\left[\frac{n+1}{2}\right]}k = \left(\left[\frac{n+1}{2}\right]+1\right)\left(\left[\frac{n+1}{2}\right]+2\right)\dots n,$$

which yields

$$2^{\left[\frac{n+1}{2}\right]}k(n+1) = \left(\left[\frac{n+1}{2}\right]+1\right)\left(\left[\frac{n+1}{2}\right]+2\right)\dots(n+1),$$

where $k(n+1)$ is odd.

Case(ii) If $n$ is odd, then

$$2^{\left[\frac{n}{2}\right]}k = \left(\left[\frac{n}{2}\right]+1\right)\left(\left[\frac{n}{2}\right]+2\right)\dots n$$

implies

$$2^{\frac{n-1}{2}}k = \left(\frac{\mathrm{n}+1}{2}\right)\left(\frac{\mathrm{n}+3}{2}\right)\dots n,$$

which yields

$$2^{\frac{n+1}{2}}k = \left(\frac{\mathrm{n}+1}{2}+1\right)\left(\frac{\mathrm{n}+1}{2}+2\right)\dots(n+1).$$

So every group of order $n(n-1)\dots(\left[\frac{n}{2}\right]+1)$ has a $Sylow\ 2-subgroup$ of order $2^{\left[\frac{n}{2}\right]}$.